\documentclass[12pt]{article}
\usepackage{amsmath,amsfonts,amsthm,amssymb}
\usepackage[T2B]{fontenc}
\usepackage[english,russian]{babel}
\usepackage[cp1251]{inputenc}
\usepackage{babel}
\usepackage{graphics}

\exhyphenpenalty=10000 \theoremstyle{plain}

\theoremstyle{plain}
\newtheorem{theorem}{Theorem}
\newtheorem{MainResult}{Theorem}
\newtheorem{Proposition}{Statement}
\newtheorem{Lem}{Lemma}
\newtheorem*{fact}{Fact}
\newtheorem*{corollary}{Corollary}
\theoremstyle{remark} \theoremstyle{definition}
\newtheorem{Def}{Definition}
\newtheorem*{MainDef}{Definition}
\newtheorem{Rem}{Remark}

\author {N.\,V.\,Gravin }
\title{Non-degenerate colorings in the Brook's Theorem}
\date{}
\begin{document}

\selectlanguage{english}  \maketitle

\newcounter{fig}
\addtocounter{fig}{1}

\abstract{

Let $c\geq 2$ and $p\geq c$ be two integers. We will call a proper
coloring of the graph $G$ a \textit{$(c,p)$-nondegenerate}, if for
any vertex of $G$ with degree at least $p$ there are at least $c$
vertices of different colors adjacent to it.

In our work we prove the following result, which generalizes Brook's
Theorem. Let $D\geq 3$ and $G$ be a graph without cliques on $D+1$
vertices and the degree of any vertex in this graph is not greater
than $D$. Then for every integer $c\geq 2$ there is a proper
$(c,p)$-nondegenerate vertex $D$-coloring of $G$, where
$p=(c^3+8c^2+19c+6)(c+1).$

During the primary proof, some interesting corollaries are derived.}

\bigskip


\section*{Introduction}

We follow the terminology and notations of the book [5] and consider
finite and loopless graphs. As in [5], $\delta(G)$ and $\Delta(G)$
denote the minimal and the maximal degree of a graph $G$
respectively. For a vertex $v\in V(G)$ the \textit{neighborhood} of
$v$ in G is $N_{G}(v)=\{u\in V(G) :$ $u$ is adjacent to $v$ in
$G$$\}$. Vertices in $N_{G}(v)$ are called neighbors of $v$. Also
$|S|$ denotes the cardinal number of a set $S$.

For an integer $k>0$, let $\overline{k}=\{1,2,\cdots,k\}$. A
\textit{proper k-coloring} of a graph $G$ is a map $c : V(G)\mapsto
\overline{k}$ such that if $u,v\in V(G)$ are adjacent vertices in
$G$, then $c(u)\neq c(v)$. Let $c$ is a proper $k$-coloring of $G$
and a set $V'\subseteq V(G)$, then by $c(V')$ we denote a
restriction of the map $c$ to the set $V'$, so we get a proper
$k$-coloring of the induced graph $G(V')$.

A proper vertex $k$-coloring is a \textit{proper conditional
$(k,c)$-coloring}, if for any vertex of degree at least $c$ there
are at least $c$ different colors in its neighborhood. This notion
for $c=2$ appeared in the works [3] and [4] as a dynamic coloring.
But results obtained there were not the Brook's Theorem
generalizations, because a number of colors in which graph was
colored is bigger then it is in the Brook's Theorem.

Further development of this theme can be found in the work\,[6]
where the definition of a conditional coloring has been given for
the first time. In this paper authors remarked that it would be
interesting to know an analogous of Brook's Theorem for conditional
colorings. But the problem of finding such an analogous seems to be
too hard in such formulation. Let us show the consideration, which
lets one to think about changing the statement. If there is a vertex
of degree $c$ in the graph, then in any $(k,c)$-coloring all its
neighbors will be colored with different colors and it means that we
can replace this vertex by $c$-hyperedge on its neighborhood.
Repeating such transformations with a graph, we can obtain any graph
with $c$-hyperedges and simple edges. So we can extend our results
of just proper colorings on such graphs. But a graph with hyperedges
is a complicated object for investigation concerning proper
colorings. Even for $c = 3$ one can easily construct a big variety
of graphs of the maximal degree $D$(for sufficiently large $D$)
which have no conditional $(D+100,3)$-coloring just by drawing the
complete graph on a $D+101$ vertices and changing some of its
triangle subgraphs to $3$-hyperedge in such a way that all vertices
will have degree not greater than $D$. So it seems to us natural to
change a little definition of the conditional coloring. The crucial
consideration, which allowed us to get serious progress in this
field, is that we demand another condition of non-degenerateness of
a proper coloring. We will call this demand the
$(c,p)$-nondegenerateness.

\begin{MainDef}
Let $c\geq 2$ and $p\geq c$ be positive integers. We call a vertex
coloring of a graph $G$ $(c,p)$\textit{-nondegenerate} if for any
vertex, with degree at least $p$, there are at least $c$ vertices of
different colors among all its neighbors.
\end{MainDef}

So, speaking informally, we impose the requirement of
nondegenerateness only to vertices of a large degree. But with such
a weaker new requirement, we can state and prove stronger and more
general theorem.

\begin{MainResult}
\label{t1}

Let $D\geq 3$ and $G$ be a graph without cliques on $D+1$ vertices
and $\Delta_{G}\leq D$. Then for every integer $c\geq 2$ there is a
proper $(c,p)$-nondegenerate vertex $D$-coloring of $G$, where
$p=(c^3+8c^2+19c+6)(c+1).$
\end{MainResult}

One of the main steps in the proof of the \textit{theorem \ref{t1}}
is the following \textit{theorem \ref{t2}}, which by itself appears
to be an interesting result.

\begin{MainResult}
\label{t2}

Let $G$ be a graph with no cliques on $D+1$ vertices with
$\Delta_{G}\leq D$. And let $D=\sum\limits_{i=1}^{c+1}\alpha_i$,
where $\alpha_i\geq 2$ are integer numbers. Then in the set $\Xi$ of
all colorings of $G$ with $c+1$ colors there is a coloring $\xi$
such that:\\
1) $\Phi(\xi)={\min\limits_{\psi\in \Xi}\Phi(\psi)}$, where
$\Phi=\sum\limits_{i=1}^{c+1}\frac{f_i}{\alpha_{i}}$ and $f_i$ is
a number of edges in $G$ connecting vertices of the $i$-th color.\\
2) For any $1\leq i\leq c+1$, there are no cliques on $\alpha_{i}+1$
vertices of the $i$-th color in $\xi$.
\end{MainResult}

In particular, there is a direct corollary from the \textit{theorem
\ref{t2}}, which is similar to the result, obtained by L.\,Lovasz in
the paper [2].

\begin{corollary}
Let $G$ be a graph with no cliques on $D+1$ vertices with
$\Delta_{G}\leq D$. And let $D=\sum\limits_{i=1}^{k}\alpha_i$, where
$\alpha_i\geq 2$ are integer numbers. Then the set $V(G)$ can be
splited into $k$ subsets $V_1$, $V_2$,... ,$V_k$ so that for any
$i\in [1,k]$ there are no cliques on $\alpha_i+1$ vertices in
$G(V_i)$ and $\Delta_{G(V_i)}\leq \alpha_i$.
\end{corollary}

\section*{Main theorem proof}

\begin{Rem}
\label{r-15}

The $(c,p)$-nondegenerateness of a coloring is a rather strong
condition even in a case of a bipartite graph (and coloring it with
$p$ colors), since it is not easy to prove a statement analogous to
the \textit{theorem \ref{t1}}. And if we want to get a
$(c,p)$-nondegenerate proper $D$-coloring of a bipartite graph but
do not bound the maximal degree of this graph, then the statement of
\textit{theorem \ref{t1}} doesn't hold for $c=2$ and every $p$.

\textit{Contrary instance:}

We take a set $S_1$ consisting of $(p-1)D+1$ elements as the first
part of $G$. As the second part of $G$, we take the set of all
$p$-element samplings from $S_1$ and join every such sampling with
all its elements in $S_1$ (see fig.\,\arabic{fig}). If we try to
color $G$ with $D$ colors, then by the Dirichlet principle in the
set $S_1$ one can find $p$ vertices of the same color and this means
that for correspondent $p$-element sampling in $S_2$ the
$(2,p)$-nondegenerate condition does not hold.

\begin{center}
\includegraphics{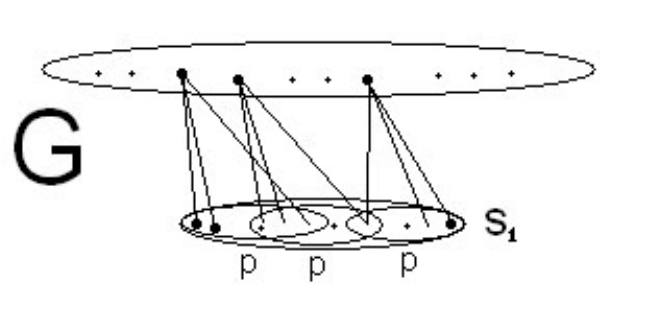}

\textit{fig. \arabic{fig}.}

\end{center}
\addtocounter{fig}{1}

\end{Rem}

\begin{Rem}
\label{r-1}

Unfortunately our estimation $p(c)=(c^3+8c^2+19c+6)(c+1)$ gives
rather large value for a small $c$. It is quite possible that using
our proof method one can get a better estimation, but it is
impossible to get an estimation asymptotically better than
$c^4(1+O(c^{-1}))$ using only our method.

\end{Rem}

\begin{theorem}
\label{t1}

Let $D\geq 3$ and $G$ be a graph without cliques on $D+1$ vertices
and $\Delta_{G}\leq D$. Then for every integer $c\geq 2$ there is a
proper $(c,p)$-nondegenerate vertex $D$-coloring of $G$, where
$p=(c^3+8c^2+19c+6)(c+1).$
\bigskip

\begin{Proposition}
\label{p0}

Without loss of generality graph $G$ may be thought of as a graph
containing no vertices of degree less than $p$.

\begin{proof}

The following operation can be done with $G$: take two copies of $G$
and join in this copies all pairs of similar vertices with degree
less than $p$ (see fig. \arabic{fig}).

\begin{center}
\includegraphics{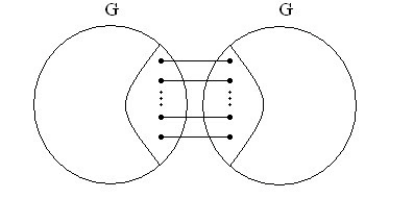}

\textit{fig. \arabic{fig}.}

\end{center}
\addtocounter{fig}{1}

Obtained graph satisfies all the conditions of \textit{theorem
\ref{t1}}. Also let us notice that if we get a $(c,p)$-nondegenerate
proper $D$-coloring of the obtained graph then we get the same for
an every copy of $G$. We repeat this operation while there is
vertices of degree less than $p$. We repeat this operation a finite
number of times because, by every execution of such operation, we
increase the smallest degree of a graph.

\end{proof}
\end{Proposition}

\bigskip

\begin{proof}

The proof of \textit{theorem \ref{t1}} consists of two parts. In the
first part we reduce our theorem to some lemma (see lemma \ref{l1}).
And in the second part we prove this lemma.

\section*{The first part.}

Choose such a number $\alpha_i$ for every $i\in\{1,2,...,c+1\}$,
that $\alpha_{i}=\left\lceil\frac{D}{c+1}\right\rceil$ or
$\alpha_{i}=\left[\frac{D}{c+1}\right]$ and
$\sum\limits_{i=1}^{c+1}\alpha_i=D$ (it is clear that we can choose
such a set of $\alpha_i$). Consider for every coloring $\xi$ with
colors $\{1,2,...,c+1\}$ a function $\Phi(\xi)$ which is determined
as follows:
$\Phi(\xi):=\sum\limits_{i=1}^{c+1}\frac{f_i}{\alpha_{i}}$, where
$f_i$ is a number of edges connecting vertices of the $i$-th color
in the coloring $\xi$. Then consider those colorings of the graph
$G$ with $c+1$ colors for which $\Phi$ reaches its minimum. Denote
such a set of colorings as $G_c$. It is obvious that $G_c$ is not
empty. Then for any coloring $\xi$ from the set $G_c$ the following
statements hold:

\begin{Proposition}
\label{p1}

For every color $i\in\{1,2,...,c+1\}$ in $\xi$ and every $i$-th
color vertex $v$ of $G$ a number of vertices adjacent to $v$ of the
$i$-th color does not exceed $\alpha_i$.

\begin{proof} Suppose the statement is false. Then from the condition that
$\sum\limits_{j=1}^{c+1}\alpha_j=D$ there can be found a color $j$
such that $v$ is adjacent in the graph $G$ to less than $\alpha_j$
$j$-th color vertices. So by recoloring $v$ with the color $j$ we
arrive at a contradiction.
\end{proof}
\end{Proposition}

\begin{Proposition}
\label{p2}

If some vertex $v$ of the $i$-th color in the coloring $\xi$ of $G$
is adjacent to exactly $\alpha_i$ vertices of the $i$-th color then
$v$ is adjacent to exactly $\alpha_j$ vertices of the $j$-th color
for every color $j$.

\begin{proof}

Assume the opposite to the \textit{statement \ref{p2}} assertion.
Then by condition that $\sum\limits_{k=1}^{c+1}\alpha_k=D$ there can
be found a color $j'\neq i$ such that $v$ is adjacent in $G$ to less
than $\alpha_{j'}$ vertices of the $j'$-th color. So by recoloring
$v$ with the color $j'$ we arrive at a contradiction.

\end{proof}
\end{Proposition}

\begin{Proposition}
\label{p3}

If the vertex $v$ of the $i$-th color in the coloring $\xi$ of the
graph $G$ is adjacent to at least one vertex of the $i$-th color
then it is adjacent to at least one vertex of any other color.

\begin{proof}
Suggesting that statement fails we arrive at a contradiction with
minimality of $\Phi(\xi)$ by recoloring $v$ with the color to which
$v$ is not adjacent.
\end{proof}
\end{Proposition}

\bigskip

We are going to prove now that there is a coloring in the coloring
set $G_c$ with no $\alpha_{i}+1$ cliques in $G$ of the $i$-th color.
We will call such cliques the \textbf{large} cliques.

\medskip

Due to the \textit{statement \ref{p1}} there can not be bigger
cliques of the $i$-th color in $G$ for any coloring from $G_c$.

For every coloring $\xi$ in $G_c$ denote as $\phi(\xi)$ a number of
\textit{large} cliques in $\xi$. Denote by $\Omega$ the set of all
colorings in $G_c$ with the smallest number of the \textit{large}
cliques. Let $\phi>0$ for all colorings in $\Omega$.

Then using the \textit{statement \ref{p2}} we get:

\begin{Proposition}
\label{p35}

If we take a vertex $v$ from some \textit{large} clique in some
coloring $g_c\in \Omega$ and recolor this vertex with any other
color then an obtained coloring $g_{c}'\in G_c$ and
$\phi(g_{c}')\leq\phi(g_{c})$.

\end{Proposition}

In \textit{statement \ref{p35}} we took $\phi$ to be the minimal on
colorings from $G_c$, so a number of \textit{large} cliques
shouldn't change. And it means that a \textit{large} clique should
appear on vertices of the color with which we recolored $v$, besides
we get $g_{c}'\in\Omega$.

\begin{Proposition}
\label{p4}

Let coloring $\xi_1\in\Omega $ and $ \phi(\xi_1)>0$. Let $C_1$ be a
\textit{large} clique of the $i$-th color. Consider the induced
subgraph $G_{ij}$ of $G$ on all vertices of the $i$-th and $j$-th
colors. Then connectivity component containing $C_1$ in the graph
$G_{ij}$ constitute a complete graph on $\alpha_{i}+\alpha_{j}+1$
vertices.

\begin{proof}

Recolor an arbitrary vertex $v_1\in C_1$ with the color $j$.
According to the \textit{statement \ref{p35}} we get a new coloring
$\xi_2\in\Omega$. And $v_1$ should get in some \textit{large} clique
$C_2$ of the $j$-th color. Recolor some distinct from $v_1$ vertex
$v_2$ in the clique $C_2$ with the color $i$. Again according to the
\textit{statement \ref{p35}} we get a new coloring $\xi_3\in\Omega$
in which $v_2$ necessarily should get in some \textit{large} clique
$C_3$ of the $i$-th color. And so on: we recolor vertices in such a
manner until we get the \textit{large} clique a part of which we
have already considered (see fig. \arabic{fig}, where four
recolorings have been done and $\alpha_{i}=\alpha_{j}=3$).

1.a) At the end we came back to a part of the clique $C_1$ and a
number of recolorings is greater than two, i.e. the last coloring is
$\xi_{k}$ where $k\geq 3$. Recolor in the coloring $\xi_1$ some
another than $v_1$ vertex $v$ in the clique $C_1$ with $j$ color.
According to the \textit{statement \ref{p35}} we get a
\textit{large} clique containing $v_k$ and $v$ of the color $j$ and
therefore the following holds: any vertex $v\in C_1$, where $v\neq
v_1$, is adjacent to all vertices in $C_{k}$ except $v_{k-1}$.

Draw the following conclusion:

Any vertex $u\in C_k$, where $u\neq v_{k-1}$, is adjacent to all
vertices in $C_1$ except $v_1$.

Recolor in $\xi_1$ vertex $v\in C_1$, $v\neq v_1$ with the $j$-th
color and then recolor some vertex $u\in C_k$ distinct from
$v_{k-1}$ and $v_k$ with the $i$-th color (we can choose such a
vertex $u$ because of $\alpha_i\geq 2$ and $\alpha_j\geq 2$). So we
get a coloring $\xi'\in G_c$ with a smaller value of $\phi$ as $u$
is adjacent to all vertices in $C_1$ except $v_1$. The following
figure \arabic{fig} is called upon to illustrate process of
recolorings for $k=4$ and $\alpha_{i}=\alpha_{j}=3$.

\begin{center}

\includegraphics{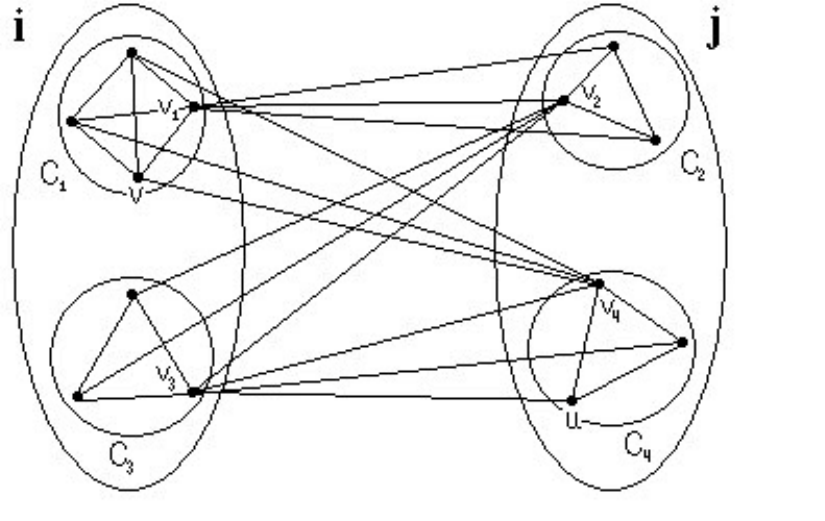}

\textit{fig. \arabic{fig}.}
\end{center}

1.b) Point out that if it was only two recolorings and we came back
to a part of the clique $C_1$ then the vertex $v_2$ is adjacent to
all vertices in $C_1$ and so by recoloring in $\xi_1$ of any vertex
in the \textit{large} clique $C_1$ with the $j$-th color we will get
by the \textit{statement \ref{p35}} a new \textit{large} clique of
the $j$-th color containing $C_2\setminus\{v_1\}$. So $G(C_{1}\cup
C_{2})$ is a complete graph. By arbitrary choice of the $v_1$ and
$v_2$ and by the fact that $G(C_{1}\cup C_{2})$ is a
$\alpha_i+\alpha_j+1$ size clique it follows that vertices of the
set $C_1\cup C_2$ are not adjacent to the rest vertices of the
$i$-th and $j$-th colors.

2) If we interrupted the process of recolorings on a clique $C_l$
where $l$ not necessary equals to $1$ then by above reasoning it is
clear that (we can assume that we start the process from $\xi_l$)
$C_{l}\cup C_{l+1}$ constitute a clique in $G$. And so we get $l=1$,
because vertices from $C_{l}\cup C_{l+1}$ and the rest vertices of
the $i$-th and $j$-th colors are not adjacent.

\end{proof}

\begin{Rem}
\label{r1}

Note that at the \textit{statement \ref{p4}} proof we make essential
use of $\alpha_i\geq 2$ and $\alpha_j\geq 2$. In other case we just
could not choose a vertex distinct from all $v_i$.

\end{Rem}
\end{Proposition}

\begin{Proposition}
\label{p5}

In any coloring $g_c\in\Omega$ there are no \textit{large} cliques.

\begin{proof}

There is a coloring $g_c\in\Omega$ with a \textit{large} clique $C$
on vertices of the $i$-th color. Without loss of generality suppose
that $i=1$. Apply the \textit{statement \ref{p4}} to the first and
the second colors. We get a complete graph containing $C$ on
$\alpha_{1}+\alpha_{2}+1$ vertices of the first and the second
colors. We can split in arbitrary way this complete graph into two
parts of the first and the second colors with correspondent sizes
$\alpha_{1}+1$ and $\alpha_2$ preserving remain coloring of the
graph and an obtained coloring would also lay in $\Omega$. By the
\textit{statement \ref{p4}} and above consideration applying to the
first and the $i$-th color ($i\in [2,c+1]$) it's easy to show the
presence of a complete subgraph of $G$ on
$1+\sum\limits_{j=1}^{c+1}\alpha_{j}$ vertices, i.e. the complete
subgraph on $D+1$ vertices -- contradiction with the condition of
\textit{theorem \ref{t1}}.

\end{proof}
\end{Proposition}

\begin{Rem}
\label{r15}

In fact we have just now proved the \textit{theorem \ref{t2}}. Also
note that desired in the \textit{theorem \ref{t1}} coloring $\xi$
assign a partition of all vertices of the graph into required in the
\textit{corollary} sets.
\end{Rem}

\begin{Rem}
\label{r2}

Consider the particular coloring $g_c\in\Omega$. We have just shown
that in $g_c$ there is no \textit{large} clique. So using the
Brook's theorem for any color in $g_c$ we can get a proper
$\alpha_i$-coloring of $i$-th color vertices, so as a result we can
get a proper coloring of $G$ with $D$ colors
($\sum\limits_{j=1}^{c+1}\alpha_{i}=D$). If a vertex in the coloring
$g_c$ is adjacent to some vertex of its color, then by
\textit{statement \ref{p3}} there should be at least $c+1$ vertices
of different colors in the neighborhood of such a vertex. In other
words the main problem we have to solve is to satisfy the condition
of $(c,p)$-nondegeneration for ``singular'' vertices, i.e. vertices
not adjacent to its and some other colors in the coloring $g_c$. In
fact, if $G$ is a bipartite graph then the theorem about $G$ proper
$(c,p)$-nondegenerate coloring with $D$ colors would be none trivial
fact. And a proof of the theorem for the case of a bipartite graph
would show you a difficulty and specificity of the problem.

\end{Rem}
\bigskip

Consider a coloring $g_{c}\in\Omega$ and consider in it all vertices
adjacent to less than $c-1$ different colors. Denote a set of all
such vertices by $\Upsilon$. Notice that every vertex $v\in\Upsilon$
has no adjacent to it vertices of the same as $v$ color in the
coloring $g_c$ and there is another color such that $v$ is not
adjacent to the vertices of this color. So we can change color of
$v\in\Upsilon$ into another such that obtained coloring as before
would be in $\Omega$. Moreover we can change color of any part of
vertices from $\Upsilon$ of an $i$-th color so that obtained
coloring will be in $\Omega$ (of course we could recolor this
vertices with different colors). For every vertex $v\in\Upsilon$
there can be found a color in $g_c$ such that $v$ is adjacent to at
least $\lceil\frac{p}{c-1}\rceil$ vertices of this color. So we can
divide $\Upsilon$ into $c+1$ sets
$\theta_1,\theta_2,...,\theta_{c+1}$, in such a way that every
vertex from $\theta_i$ is adjacent to at least
$\lceil\frac{p}{c-1}\rceil$ vertices of the $i$-th color.

\bigskip

Denote by $H_i$ for all $i\in [1,c+1]$ the induced subgraph of $G$
on the vertices of the $i$-th color in the coloring $g_c$.

\begin{Proposition}
\label{p6}

For any vertex $v\in H_i$ the following inequality holds:
$$ \lceil d_{H_{i}}(v)+\frac{d_{G(\theta_{i}\cup\{v\})}(v)}{c+2} \rceil\leq\alpha_{i}.$$

\begin{proof}
Consider a set $E_v$ of all edges in the graph $G$ with one end at
$v$. It's obviously that $|E_v|\leq D$. Consider a set $E_1$ of all
edges from $E_v$ which has the second end vertex distinct from $v$
not laying in $\theta_i$. Let from $v$ there lead less than
$\frac{\alpha_j}{\alpha_i}d_{H_{i}}(v)$ edges of the set $E_1$ to a
color $j$ distinct from $i$. Then we change the color of all
vertices of the $j$-th color of the set $\theta_i\subseteq \Upsilon$
in such a way that an obtained coloring will be in $\Omega$. Clearly
we recolored these vertices not with the color of $v$, so
$d_{H_{i}}(v)$ doesn't change in the obtained coloring. If we
recolor $v$ in the new coloring with the $j$-th color then a
magnitude less than
$\frac{\frac{\alpha_j}{\alpha_i}d_{H_{i}}(v)}{\alpha_j}-\frac{d_{H_{i}}(v)}{\alpha_i}=0$
is added to the value of $\Phi$, thus we find a coloring with a
smaller value of $\Phi$ and so we arrive at a contradiction.

\begin{center}
\addtocounter{fig}{1}
\includegraphics{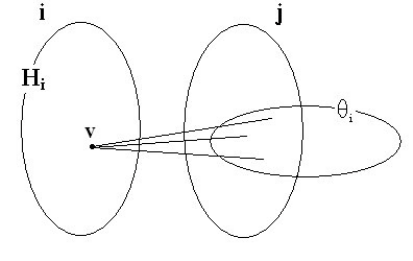}

\textit{fig. \arabic{fig}.}
\end{center}

So we can get the following lower bound on the number of edges
coming from $v$:

$|E_v|\geq d_{H_{i}}(v)+d_{G(\theta_{i}\cup\{v\})}(v)+
\sum\limits_{j\neq
i}\frac{\alpha_j}{\alpha_i}d_{H_{i}}(v)=\sum\limits_{j=1}^{c+1}\frac{\alpha_j}{\alpha_i}d_{H_{i}}(v)+d_{G(\theta_{i}\cup\{v\})}(v)=
\frac{D}{\alpha_i}d_{H_{i}}(v)+d_{G(\theta_{i}\cup\{v\})}(v)~.$

By definition $|E_{v}|\leq D$. So we get:

$$D\geq \frac{D}{\alpha_i}d_{H_{i}}(v)+d_{G(\theta_{i}\cup\{v\})}(v)\Rightarrow \alpha_i\geq d_{H_{i}}(v)+\frac{\alpha_i}{D}d_{G(\theta_{i}\cup\{v\})}(v)~.$$

Then by using the fact that $\alpha_i \geq [\frac{D}{c+1}]$ and
$D\geq (c^3+8c^2+19c+6)(c+1)$ we get
$\frac{\alpha_i}{D}>\frac{1}{c+2}$. So we get $$\alpha_i\geq
d_{H_{i}}(v)+\frac{d_{G(\theta_{i}\cup\{v\})}(v)}{c+2}~.$$
\end{proof}

\end{Proposition}

\section*{The second part.}

\begin{Lem}
\label{l1} Let there are given two non empty sets $A$ and $B$ and a
connected graph $H=(A\cup B, E)$. And let $\cal{G}$ denotes the
induced subgraph $H(B)$. Define $d_{_A}(v), v\in B$ to be a number
of edges coming from $v$ to the set $A$. Let the graph $H$ satisfy
the following conditions:

1) every two vertices of $A$ are not joint with an edge;

2) the degree of every vertex from $A$ in the graph $H$ is at least
$d$, where $d=q^3+2q^2-q-8$ and $q\geq 4$;

3) for any vertex $v\in B$, the following inequality holds:
$$ d_{\cal{_G}}(v)+ \lceil \frac{d_{_A}(v)}{q}\rceil \leq d. \eqno(1)$$

Then the graph $\cal{G}$ could be properly colored with $d$ colors
in such a way, that for any vertex $v\in A$ among all its neighbors
in $B$ there are vertices of at least $q$ different colors.

\begin{Rem}
\label{r36}

$(c+2)^3+2(c+2)^2-(c+2)-8=c^3+8c^2+19c+6$.

\end{Rem}

\begin{Rem}
\label{r35}

In the \textit{lemma \ref{l1}}, the set $B$ denotes $H_i$ from the
\textit{first part}, the set of vertices $A$ denotes $\theta_i$ from
the \textit{first part}. Also it makes no difference for us whether
there are any edges between vertices in $\theta_i$. We only need to
know to which vertices in $H_i$ vertices in $\theta_i$ are adjacent
to, because we will color vertices only in $H_i$.

As $q$ in \textit{lemma \ref{l1}}, we denoted the value of $c+2$
from the \textit{first part} and as $d$ we denoted the value of
$\alpha_i$. Via $H$ in the \textit{lemma \ref{l1}} we denoted the
graph $G(\theta_{i}\cup H_{i})-E(G(\theta_i))$. By definition of the
set $\theta_i$ from any vertex $v\in\theta_i$ there comes at least
${\frac{p}{c-1}>q^3+2q^2-q-8}$ edges to the set $V(H_{i})$.

We suppose in the \textit{lemma \ref{l1}} that the graph $H$ is
connected (in other case it is sufficient to prove the lemma's
statement for every connectivity component). Furthermore we can
assume that $\theta_i$ is not empty, otherwise we have just to prove
the Brook's Theorem because of we need to color properly graph $H_i$
with $\alpha_i$ colors, and we know that in $H_i$ there are no
complete subgraphs on $\alpha_{i}+1$ vertices ( in $H_i$ there are
no \textit{large} cliques) and $d_{H_{i}}=\lceil
d_{H_{i}}(v)+\frac{d_{G(\theta_{i}\cup\{v\})}(v)}{c+2}
\rceil\leq\alpha_{i}$. Thus, all the conditions of \textit{lemma
\ref{l1}} are satisfied for the sets $B=V(H_{i})$ and $A=\theta_i$.

Suppose the \textit{lemma \ref{l1}} has been already proven. Then,
if we color for every $i$ the subgraph $H_i$ in the coloring $g_c$
of $G$ in a proper way with a new $\alpha_i$ colors such that every
vertex from $\theta_i$ would be adjacent to vertices of at least $c$
different colors then we get a proper $D$-coloring of the whole
graph $G$. At that time the vertices from the set
$\Upsilon=\bigcup\limits_{i=1}^{c+1}\theta_i$ would be adjacent to,
at least, $c$ vertices of different colors. Moreover in accordance
with the definition of $\Upsilon$ all the vertices from the set
$V(G)\setminus\Upsilon$ would be adjacent to at least $c$ vertices
of different colors. Thus, we reduce the \textit{theorem \ref{t1}}
to the \textit{lemma \ref{l1}}.

\end{Rem}

\begin{Rem}
\label{r3}

The second part is devoted to the proof of \textit{lemma \ref{l1}}.
So to avoid a misunderstanding for a coincidence of notations let us
say that notations from the \textbf{first part} have no connection
with notations from the \textbf{second part}.

\end{Rem}

\begin{Rem}
\label{r355}

In the assertion of the \textit{lemma \ref{l1}} it is possible to
change $q$ to $q-2$, but we will not do this for the sake of
calculation convenience.
\end{Rem}

\bigskip

\begin{proof}[Proof of the Lemma \ref{l1}]

Suppose that assertion of the \textit{lemma \ref{l1}} fails. Then,
consider the smallest for a number of vertices graph for which all
the assumptions of the \textit{lemma \ref{l1}} holds but the
statement of the \textit{lemma \ref{l1}} fails.

\begin{Def}
We will call a \textit{permissible} the set $S_i\subseteq B$ if
$S_i\subset N_{_H}(v_i)$, where $v_{i}\in A$ and $|S_i|=q$. A set of
all samplings of \textit{permissible} sets for all
$i\in\{1,2,...,|A|\}$ we will denote by $\Lambda$.

\end{Def}

\medskip

The assertion of our lemma abides by the following fact:

\begin{fact}
For every vertex $v_i$ in $A$ we can choose a \textit{permissible}
set $S_i$ in such a way that if we add to the edges set $E(\cal{G})$
all complete graphs on sets $S_i$ where $i\in\{1,2,...,|A|\}$ then
it is possible to color vertices of the obtained graph
$\widetilde{\cal{G}}$ properly with $d$ colors.

\end{fact}

\begin{Rem}
\label{r4}

We will consider $\widetilde{\cal{G}}$ as a graph with multiedges.
\end{Rem}
\begin{Rem}
\label{r5}

So we get an equivalent statement of the \textit{lemma~\ref{l1}}.
\end{Rem}

\begin{Rem}
\label{r6}

In the new formula, it is convenient to make some reduction with a
graph as follows:

Let there be a vertex $\hat{v}$ of degree $d$ in a graph
$\widetilde{\cal{G}}$, then it is possible to ``delete'' this vertex
from the graph $\widetilde{\cal{G}}$ and prove a statement of the
\textit{fact} for the graph $\widetilde{\cal{G}}\setminus\hat{v}$.
\end{Rem}

\begin{Def}
We will say that $\hat{v}$ is recursively deleted from
$\widetilde{\cal{G}}$ if there is a sequence of reductions described
above with the last $\hat{v}$ reduction. We will call a graph
$\widetilde{\cal{G}}$ to be a recursive one, if it reduces to the
empty graph.
\end{Def}

\begin{Rem}
\label{r7}

Let us explain why we call such a reduction as a recursion. The
matter is that if a graph reduces to the empty one then we will
color it just by recursion.
\end{Rem}

Actually we will prove the following stronger fact:

\textit{Instead of the statement that $\widetilde{\cal{G}}$ is
properly colored with $d$ colors, we will prove that
$\widetilde{\cal{G}}$ is a recursive with respect to coloring it
with $d$ colors.}

\medskip

Return to the lemma's proof and more specifically to the proof of
the stronger \textit{fact}. Denote as $S$ the set of vertices from
$B$ which are adjacent to at least one vertex in $A$.

\medskip

Prove that for the graph $H$ the strengthened \textit{fact} holds in
assumption that $H$ is the minimal for number of vertices graph for
which the statement of the \textit{lemma \ref{l1}} fails. Thus, we
will arrive at a contradiction and so we will prove the
\textit{lemma \ref{l1}}.

\bigskip

\begin{Def}

Define for any vertex $v$ from the set $B$ the magnitude $$ L(v)
:=d_{\cal{_G}}(v)+ \frac{d_{_A}(v)}{q+1}~.$$

\end{Def}

\begin{Rem}
\label{r8}

Notice that if we choose a sampling $S_i$ at random (independently
for any vertex $v_i$ where all possible variants of the set $S_i$
are equiprobable), then the distribution average of a variate of the
degree in the graph $\widetilde{\cal{G}}$ for any vertex from the
set $S$ is not greater than
$d_{\cal{_G}}(v)+d_{_A}(v)(q-1)\frac{q}{q^3+2q^2-q-8}$, i.e. the
degree is not greater than $L(v)$ (since $q\geq 4$ then
$q^3+2q^2-q-8>q(q^2-1)$) and by the third condition of \textit{lemma
\ref{l1}} would be less than $d$. Thus, at the average the degree of
every vertex in $\widetilde{\cal{G}}$ is less than $d$. And this
gives us hope that the graph $\widetilde{\cal{G}}$ turns out to be a
recursive one, i.e. if we successively delete vertices from the
$\widetilde{\cal{G}}$ with degree less than $d$ then we arrive to
the empty graph.

\end{Rem}

For a lemma's proof completion, we only need to choose successfully
a sampling of $S_i$, i.e. to choose it in such a way that
$\widetilde{\cal{G}}$ become a recursive graph.

\begin{Def}
By the \textit{change} of some \textit{permissible} sets $S_{i_1}$,
$S_{i_2}$, ...,$S_{i_z}$ in a sampling $\lambda\in\Lambda$ to some
other \textit{permissible} sets $S'_{i_1}$,~...,~$S'_{i_z}$ we
denote a substitution of $\lambda$ for a $\lambda'\in\Lambda$, where
$\lambda'$ differs from $\lambda$ only by that the all
\textit{permissible} sets $S_{i_1}$,~...,~$S_{i_z}$ in $\lambda$ are
substituted by the other permissible sets
$S'_{i_1}$,~...,~$S'_{i_z}$. The sets $S'_{i_1}$,..., $S'_{i_z}$  we
will call the \textit{result of the change} of sets $S_{i_1}$,
$S_{i_2}$,~...,~$S_{i_z}$.

\end{Def}

Denote as $R$ the set $B\setminus S$. The degree of any vertex in
$R$ may be thought of as $d$ because by the condition of
\textit{lemma \ref{l1}} the degree of any vertex of $B$ in the graph
$\cal{G}$ is less or equal than $d$ and if degree of a vertex is
less than $d$, then it is possible to delete recursively this vertex
in $\widetilde{\cal{G}}$ for any \textit{permissible} sampling.

\begin{Proposition}
\label{p7}

Let there be given a graph $F$ such that $V(F)=S_1\cup S_2$ and
$S_1\cap S_2={\o}$, the degree of any vertex of $S_2$ in the graph
$F$ is less or equal than $D$ and in $F$ there is such a vertex
$v\in S_1$ that the graph $F(S_2\cup\{v\})$ is connected, $d_{F
}(v)< D$ and the vertex $v$ is adjacent to all the other vertices in
$S_1$. Let the graph $F(S_1)$ be properly colored with $D$ colors.
Then it is possible to extend such a vertex coloring of $F(S_1)$ to
the proper $D$-coloring of $F$.

\begin{center}
\addtocounter{fig}{1}
\includegraphics{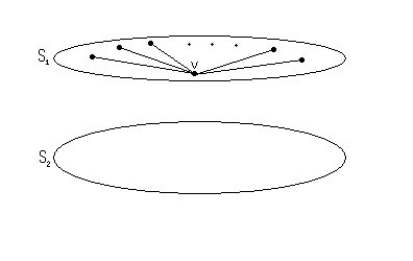}

\textit{fig. \arabic{fig}.}
\end{center}

\begin{proof}

Throw out from the graph $F$ the vertex $v$, then we get a new graph
$F'$. The set $S_1\setminus\{v\}$ has already been properly colored
with $D$ colors. One by one we recursively color properly with $D$
colors all the vertices in $S_2$, since $F(S_2\cup \{v\})$ is a
connected graph and the degree in the graph $F$ of any vertex in
$S_2$ is less or equal than $D$. Carry the obtained proper
$D$-coloring of $F'$ to $F$ and then color $v$ with some color
distinct from all the colors of vertices in $N_{F }(v)$ (it is
possible to do so since $d_{F }(v)< D$), as a result we get a proper
$D$-coloring of the graph $F$, but at that time we could probably
change the initial color of vertex $v$ in the given coloring of
$S_1$. Let us notice that all vertices in the set
$S_1\setminus\{v\}$ are colored with the colors different from the
color of $v$ in the initial coloring of $S_1$, as initial coloring
of $S_1$ was proper for the graph $F(S_1)$ and vertex $v$ is
adjacent to all the other vertices in $S_1$, moreover all the colors
of vertices in $S_1\setminus\{v\}$ differ from the color of $v$ in
the obtained proper $D$-coloring of $F$. And now if the vertex $v$
changed its color in the obtained coloring in comparison with the
given coloring of $S_1$ then we trade places of the current color of
$v$ with the color of $v$ in the initial coloring. Thus, we get a
proper $D$-coloring of $F$, but now equal on the set $S_1$ to the
initial coloring.

\end{proof}
\end{Proposition}

\begin{Def}
By the \textit{regular change} of the sets $S_i$ of a sampling
$\nu\in\Lambda$ with respect to a set $S'$, we will call such a
\textit{change} of the sets $S_i$, where $i\in [1,|A|]$, to the sets
$S'_i$, $i\in\ [1,|A|]$, that for all $i\in [1,|A|]$ the set
$S_i\cap S'$ contains the set $S'_i\cap S'$. If there exists $i\in
[1,|A|]$ such that $|S_i\cap S'|$ greater than $|S'_i\cap S'|$ then
such a \textit{regular change} we will call the
\textit{non-degenerate change}.

\end{Def}

\begin{Rem}
A \textit{Regular change} with respect to some set is a
\textit{regular change} with respect to any subset of this set, but
at that time the \textit{non-degeneracy} not necessarily preserves.
\end{Rem}

\begin{Proposition}
\label{p8}

Let there is a sampling of permissible sets
$\eta=\{S_1,S_2,...,S_{|A|}\}$ of the graph $H$ --- the smallest for
the number of vertices graph which is contrary instance for the
lemma \ref{l1} and let there are such sets $S'\subseteq S,
R'\subseteq R$ that the all vertices in $B\setminus(S'\cup R')$ are
recursively deleted from the graph $\widetilde{\cal{G}}$, for all
$u\in R ~~ d_{\widetilde{\cal{G}}}(u) =d$ and for all $u\in R' ~~
d_{\widetilde{\cal{G}}(S'\cup R')}(u)=d$.

Let $\widetilde{H'}:=\widetilde{\cal{G}}(S'\cup R')$ and $$
\sum\limits_{u\in \widetilde{H'}}d_{_{\widetilde{H'}}}(u)>
\sum\limits_{u\in \widetilde{H'}}L(u). \eqno(2)$$

Then it is possible to make a \textit{regular non-degenerate change}
of sets $S_i$ with respect to the set $S'\cup R'$ so that all the
set $B\setminus(S'\cup R')$ as before could be recursively deleted
out the graph $\widetilde{\cal{G}'}$ obtained from
$\widetilde{\cal{G}}$ as a result of this change.

\begin{proof}
We will prove this statement by induction on the set $B\setminus
(S'\cup R')$ size.

\textbf{The basis:} the case when $|B\setminus (S'\cup R')|=0$
obviously could not take place since by virtue of \textit{remark
\ref{r8}} the condition (2) doesn't hold.

\textbf{The inductive step:} let the statement holds for all numbers
less than $k$, then let us prove that it holds for the $k$.

Let $Z:=B\setminus(S'\cup R')$.

\begin{center}
\addtocounter{fig}{1}
\includegraphics{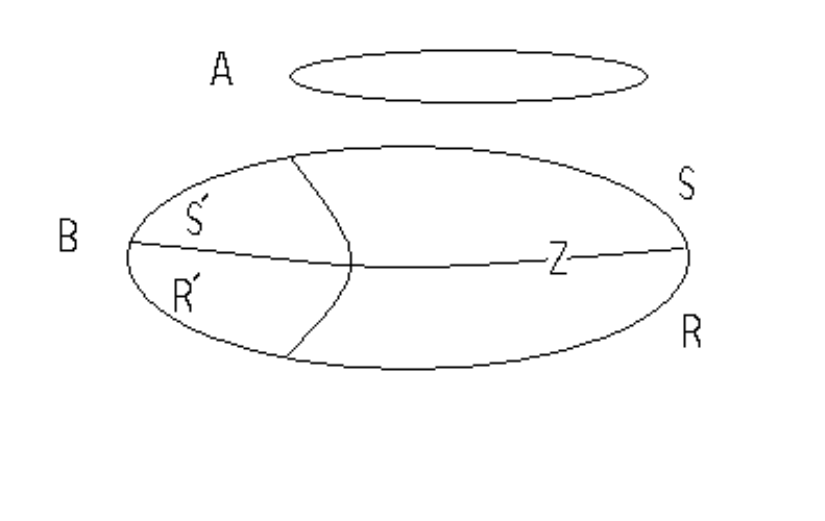}

\textit{fig. \arabic{fig}.}
\end{center}

Consider those sets $S'$ and $R'$ such that $|Z|=k$  and the
assertion of the statement fails.

\addtocounter{fig}{1}

Let us show that there is a vertex $v_i\in A$ and correspondent to
it the set $S_i$ such that it is possible to make a \textit{regular
non-degenerate change} of $S_i$ in relation to $S'$. If it is false
then for any $v_j\in A$ and correspondent to it the set $S_j$ only
two possibilities can occurred:

1) the set $S_i\cap S' = \emptyset$ (see fig. \arabic{fig});

\addtocounter{fig}{1}

2) the set $N_{\cal{_G}}(v_{j})\setminus S'\subseteq S_i$ (see fig.
\arabic{fig}).

\addtocounter{fig}{-1}

\begin{center}

\includegraphics{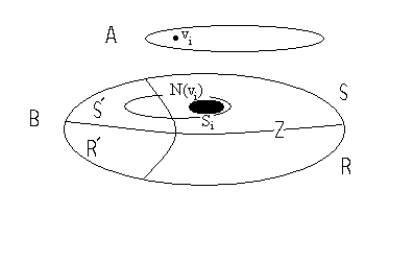}
\includegraphics{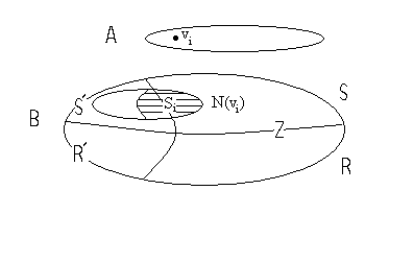}

\textit{fig. \arabic{fig}.~~~~~~~~~~~~~~~~~~~~~~~~~~~~~~~~~~~}
\addtocounter{fig}{1} \textit{fig. \arabic{fig}.}
\end{center}

In both of these cases the number of edges added to the graph
$\cal{G}$ with two ends in $S'$ reaches its minimum. Thus, for every
vertex $v\in S'\cup R'$ the following chain of inequalities take
place: $d_{\widetilde{H'}}(v)\leq E(d_{\widetilde{H'}}(v))\leq
E(d_{\widetilde{\cal{G}}}(v))\leq L(v)$, where by the $E(\cdot )$,
we denote the average of distribution of a variate with the
distribution specified in the \textit{remark \ref{r8}}. We know from
the condition (2) that $\sum\limits_{u\in
\widetilde{H'}}d_{\widetilde{H'}}(u)>\sum\limits_{u\in\widetilde{H'}}L(u)$.
So by a substitution of the inequality $d_{\widetilde{H'}}(v)\leq
L(v)$ in the previous inequality we get $\sum\limits_{u\in
\widetilde{H'}}d_{\widetilde{H'}}(u)>\sum\limits_{u\in
\widetilde{H'}}d_{\widetilde{H'}}(u)$ --- a contradiction.

Hence, there is such a vertex $v_i\in A$, that a part of its
neighborhood is contained in $Z$ but the set $S_i\cap S'\neq{\o}$
and $S_i$ does not contain this part.

Consequently, we can consider such a vertex $v\in N_{\cal{G}}(v_i)$,
that it does not lay neither in the set $S'$ nor in the set $S_i$,
but some nonempty part of $S_i$ is contained in the set $S'$. We
know that $Z$ can be recursively deleted from $\widetilde{\cal{G}}$,
so begin to \textit{recursively delete} vertices from $Z$, but do it
while it is possible to delete vertex distinct from $v$. At some
moment we should stop this process. It means that we could not
delete vertex except $v$ and so we have only $v, u_1,u_2,...,u_l \in
S, w_1,w_2,...,w_m \in R$ vertices remained in $Z$.

\medskip

Denote by $P$ the set of all remaining vertices in $Z$, and denote
by $\widetilde{I}$ induced subgraph $\widetilde{\cal{G}}(S'\cup
R'\cup P)$ of $\cal{G}$.

\medskip

Let us notice that the degree in the graph $\widetilde{I}$ for any
$u_k$ vertex, where $k\in [1,l]$, or for any $w_j$, where $j\in
[1,m]$, is at least $d$.

\medskip

Let us notice also that the degree of $v$ in $\widetilde{I}$ is less
than $d$.

\medskip

If the degree of $v$ is less than $d-q+1$ in $\widetilde{I}$, then
let us make a change of $S_i$ to a set $S'_i$ in the following way:
we take a vertex $x$ in $S_i$ which also is contained in the set
$S'\cap S_i$ (those vertex necessarily turns up as $S'\cap
S_i\neq{\o}$), then $S'_i:=\{(S_{i}\setminus\{x\})\cup \{v\}\}$, the
remaining sets of the sampling $\eta$ we do not change. Let us
notice that the change described above is a \textit{regular} and
\textit{non-degenerate} one in regard to $S'$ also it is clear that
set $Z$ will be \textit{recursively deleted} in the obtained graph
(it is clear that we can recursively delete as earlier vertices from
$Z\setminus P$ then we can recursively delete $v$, as it has degree
less than $d$, because before the change it has degree less than
$d-q+1$ and after the change the degree became not greater than
$d-1$, and then we can recursively delete all remaining vertices
from $Z$, since $Z$ has been recursively deleted from
$\widetilde{\cal{G}}$ and we drew no new edges in the graph
$\widetilde{I}(V(\widetilde{I})\setminus\{v\})$). So in this case we
have proved an inductive step.

Thus we get that the degree of $v$ is less than $d$ but at least
$d-q+1$ in $\widetilde{I}$.

\medskip

Let us prove that for the graph $\widetilde{I}$ the following
condition holds:
$$ \sum\limits_{u\in \widetilde{I}}d_{_{\widetilde{I}}}(u)> \sum\limits_{u\in \widetilde{I}}L(u). \eqno(2')$$

With the proof, we can make use of an induction assumption for the
sets $S'_{0}$ and $R'_{0}$, where
$S'_{0}:=(S'\cup{R'}\cup{P})\cap{S}$ and $R'_{0}:=(S'\cup R'\cup
P)\cap{R}$, i.e. we can make a \textit{regular non-degenerate
change} of $\eta$ in regard to $S'\cup R'\cup P$ in such a way that
the set $B\setminus(S'\cup R'\cup P)$ will be recursively deleted in
obtained graph. If a sampling had \textit{regularly changed} in
relation to $S'\cup R'\cup P$, then it is \textit{regularly changed}
in regard to $S'\cup R'$, also a composition of \textit{regular
changes} in regard to some a set is also the \textit{regular change}
in regard to this very set. Besides let us notice that in the graph
obtained by this \textit{change} all vertices from the set
$B\setminus(S'\cup R')$ will be \textit{recursively deleted}, as we
can \textit{recursively delete} at first the all vertices from
$B\setminus(S'\cup R'\cup P)$ and then we can \textit{recursively
delete} as before all vertices from $P$ since by the \textit{change}
we do not add new edges to $P$.

So we will do such \textit{changes} until either $S_i$ will be
\textit{regularly changed} in \textit{non-degenerate} way in regard
to $S'\cup R'$, or the degree of $v$ in the graph
$\widetilde{\cal{G}}$ will become less than $d-q+1$, or the degree
of any vertex from $P\setminus\{v\}$ will become less than $d$. In
the last case we can \textit{recursively delete} some more vertices
from $Z$ and for the smaller graph $\widetilde{I}$ apply the same
arguments. Here, it needs to be emphasized that some time or other
we necessarily arrive at one of this cases else we will do an
infinite number of \textit{non-degenerate regular changes} in regard
to the set $S'\cup R'\cup P$ and, hence, we will infinitely decrease
a value of the sum $\sum\limits_{i=1}^{|A|} |S_i\cap (S'\cup R'\cup
P)|$.

\begin{center}
\addtocounter{fig}{1}
\includegraphics{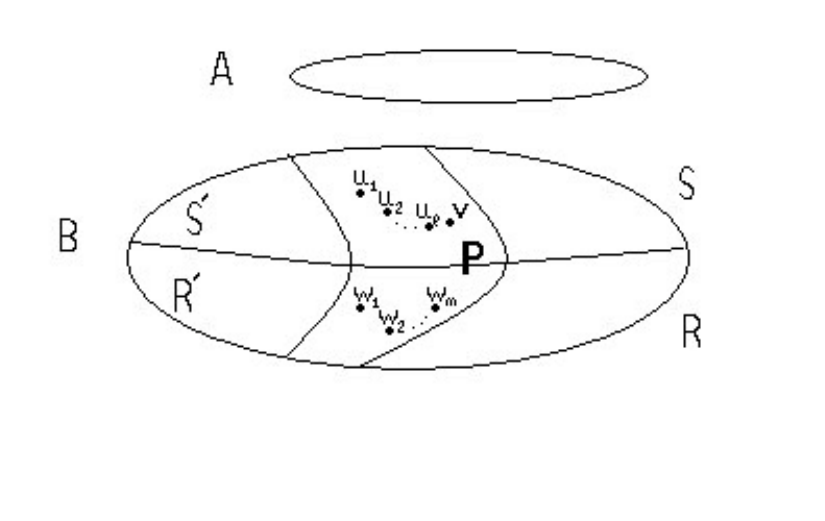}

\textit{fig. \arabic{fig}.}
\end{center}

Denote by $l'$ the number of edges coming to the vertex set $S'$
from $P$ in the graph $\widetilde{\cal{G}}$.

\medskip

By the conditions of \textit{statement \ref{p8}} that for all $u\in
R'$ $d_{\widetilde{\cal{G}}(S'\cup R')}(u)=d$ and for all $u\in R$
$d_{\widetilde{\cal{G}}}(u)=d$, there are no edges between $P$ and
$R'$.

\medskip

So to end the proof of \textit{statement \ref{p8}} we only need to
prove, that for the graph $\widetilde{I}=\widetilde{\cal{G}}(S'\cup
R'\cup P)$ the inequality (2') holds. Assume the contrary, then

$$\sum\limits_{u\in \widetilde{I}}L(u)\geq\sum\limits_{u\in \widetilde{I}}d_{_{\widetilde{I}}}(u)\geq\sum\limits_{u\in \widetilde{H'}}d_{_{\widetilde{H'}}}(u)+l'+\sum\limits_{u\in P}d_{_{\widetilde{I}}}(u)>\sum\limits_{u\in \widetilde{H'}}L(u)+\sum\limits_{u\in P}d_{_{\widetilde{I}}}(u)+l'~.$$

So we get the following:

$$\sum\limits_{u\in P}L(u)>l'+\sum\limits_{u\in P}d_{_{\widetilde{I}}}(u)'~.$$

Hence, we get the inequality:

$$l'+\sum\limits_{u\in P}d_{_{\widetilde{I}}}(u)-L(u)<0~. \eqno(3)$$

\medskip
\textbf{Let us bound the magnitude $d_{\widetilde{I}}(u_i)-L(u_i)$
for all $i\in[1,l]$}.

By definition of $L(u)$ and by virtue of $d_{\widetilde{I}}(u_i)\geq
d$ we get that for all $i\in[1,l]$ the following inequality holds:
$d_{\widetilde{I}}(u_i)-L(u_i)\geq d-d_{\cal{_G}}(u_i)-
\frac{d_{_A}(u_i)}{q+1}$. Using the inequality (1) we get:

$$d_{\widetilde{I}}(u_i)-L(u_i)\geq d-d_{\cal{_G}}(u_i)-
\frac{d_{_A}(u_i)}{q+1}\geq d_{\cal{_G}}(u_i)+\lceil
\frac{d_{_A}(u_i)}{q}\rceil-d_{\cal{_G}}(u_i)-
\frac{d_{_A}(u_i)}{q+1}~.\eqno(*)$$

Thus $d_{\widetilde{I}}(u_i)-L(u_i)\geq \lceil
\frac{d_{_A}(u_i)}{q}\rceil- \frac{d_{_A}(u_i)}{q+1}~$.

Also $d_{A}(u_i)>0$ for all $i\in[1,l]$, as $u_i\in S$ for all
$i\in[1,l]$. Let us consider two following cases:

a) $0<d_{A}(u_i)\leq q$;

b) $d_{A}(u_i)\geq q+1$.

In both of this cases the following inequality holds:

$$d_{\widetilde{I}}(u_i)-L(u_i)\geq \lceil
\frac{d_{_A}(u_i)}{q}\rceil- \frac{d_{_A}(u_i)}{q+1}\geq
\frac{1}{q+1}~.\eqno(4)$$

Let $q_1:=d_{\widetilde{I}}(v)-d+q$ then, as we have just showed it
above, $q_1>0$. Let us notice that for the vertex $v$ analogously to
calculations $(*\,)$ we can get the following inequality:
 $$d_{\widetilde{I}}(v)-L(v)\geq
q_1-q+\frac{1}{q+1}~.\eqno(5)$$

Since $w_i\in R$, where $i\in[1,m]$, $d_{\cal{G}}(w_i)=L(w_i)$,
moreover we can not \textit{recursively delete} any vertex from the
set $P\cap R$ in the graph $\widetilde{I}$. In addition using the
\textit{statement \ref{p8}} condition, that for any vertex $u\in R$
$d_{\cal{G}}(u)=d$, we get
$d_{\widetilde{I}}(w_j)=d_{\cal{G}}(w_j)=d$. And so for all $i\in
[1,m]$ the we have

$$d=d_{\widetilde{I}}(w_i)=L(w_i). \eqno(6)$$

\medskip

It now follows from (4), (5), (6), (3) that:

$$l'+l\frac{1}{q+1}-q+q_1+\frac{1}{q+1}<0~.$$

Recall now that $l$ is a number of vertices in the set $(P\cap
S)\setminus\{v\}$, i.e. the number of $u_i$. We know that $q_1\geq
1$. Then $l'(q+1)+l<(q-1)(q+1)-1$, i.e. $$(q+1)l'+l\leq q^2-3~.
\eqno(7)$$

From the inequality (7) we get two inequalities $$l\leq q^2-3
\eqno(8)$$ and $$l'\leq q-2~. \eqno(9)$$

\bigskip

\addtocounter{fig}{1}

Denote by $b_j$ $($see fig.\,\arabic{fig}$)$, where $j\in [1,r]$,
the all vertices from the set $R\cap V(\widetilde{I})$, which are
adjacent to $v$ ($r$ can be equal to $0$). Let us consider some
cases.

\bigskip

\textbf{1) $r\geq q^2-3$}.

By $C_v$ we denote the union of all connectivity components of the
graph $\widetilde{\cal{G}}(R)$, which is minimal and contains all
the vertices $b_j$, where $j\in [1,2,...,r]$. As we remark earlier,
between sets $P$ and $R'$ there are no edges, so $C_v\subseteq
R\setminus R'$. By equality (6) we have $d_{\widetilde{I}}(w_j)=d$,
where $j\in [1,m]$. Thus, vertices from the set $P\cap R$ and from
the set $Z\setminus P$ are not adjacent, and so $C_v\subseteq
\{w_1,w_2,...,w_m\}$.

\medskip

Consider, in the vertex set $S$ of the graph $\widetilde{\cal{G}}$
all adjacent to $C_v$ vertices and denote it by $W$. It is clear by
virtue of $d_{\widetilde{I}}(w_i)=d$ and $d_{\cal{G}}(w_i)=d$ that,
firstly $W\subseteq V(\widetilde{I})$, secondly $v\in W$, and
thirdly for all vertices $u\in C_v$ the equality
$d_{\widetilde{\cal{G}}(W\cup C_v)}(u)=d$ holds.

\textbf{1.1) $|W|\geq q^2-1$}.

Then $|W\cap S'|=|W|-|W\cap P\cap
 S|\geq q^2-1-(l+1)$. It is clear that $l'\geq |W\cap S'|$. Thus $l'\geq
 q^2-2-l$, i.e. $l'+l\geq q^2-2$. So we arrive at a contradiction
with inequality (7).

Thus $|W|\leq q^2-2$.
\bigskip

\textbf{ 1.2) $|W|\leq q^2-2$. }

\begin{center}

\includegraphics{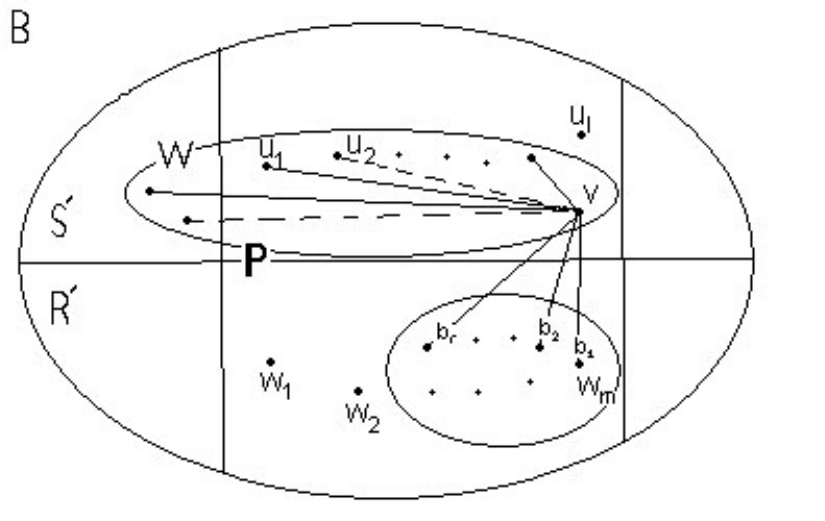}

\textit{fig. \arabic{fig}.}
\end{center}

Let us draw in the graph $H$ all the edges of the type $(u,v)$,
where $u\in W$ and $(u,v)\notin E(H)$, denote by $\Theta$ the
obtained graph. Let us verify all conditions of the \textit{lemma
\ref{l1}} for the graph $\hat{H}:=\Theta(V(H)\setminus C_{v})$ and
sets $\hat{A}:=A$, $\hat{B}:=B\setminus C_v$ and value $\hat{d}:=d$.

\textit{Condition 1)} is clear as $A$ became the same.

\textit{Condition 2)} is clear, since any vertex from $C_v$ are not
adjacent in $H$ to any vertex from the set $A$.

\textit{Condition 3)}. It is sufficient to verify inequality (1)
only for vertices from the set $B$, from which we draw any new
edges, in other words we need to verify (1) only for the set $W$. By
definition $\hat{\cal{G}}:=\hat{H}(\hat{B})$. For any vertex $u\in
W$, $u\neq v$ we added not more than one edge with the end at $u$
but also we deleted at least one edge coming from $u$ to the set
$C_v$ ($W$ is by definition the set of all vertices from $S$, which
are adjacent with at least one vertex in $C_v$). Thus (1) remains
true for all $u\in W$, $u\neq v$. The inequality (1) for $v$ holds,
as $|W|\leq q^2-2$, so we drew from $v$ not greater than $q^2-3$
edges. On the other hand the set $C_v$ by definition contains all
the $b_j$, where $j\in [1,r]$, $r\geq q^2-3$, so we deleted at least
$q^2-3$ edges with the end in $v$.

As we suppose $H$ to be a minimal by the number of vertices graph
for which the \textit{lemma \ref{l1}} doesn't hold, then
\textit{lemma \ref{l1}} holds for $\hat{H}$ which has the smaller
number of vertices. Then we can color properly the graph
$\hat{\cal{G}}$ with $d$ colors in such a way that for any vertex
$u\in A$ among its neighbors in $\hat{B}$ there would be at least
$q$ different colors. Denote by $\xi$ such a proper $d$-coloring. It
is clear that all assumptions of \textit{statement \ref{p7}} are
satisfied for the graph $\Phi:=\Theta(W\cup C_{v})$, sets $S_1:=W$
and $S_2:=C_v$ and vertex $v$. Consider a proper $d$-coloring
$\xi(W)$ of the graph $\Theta(W)$. By the \textit{statement
\ref{p7}} we can extend $\xi(W)$ to a proper $d$-coloring $\zeta$ of
the graph $\Theta (W\cup C_v)$. Let us notice that there are no
edges in the graph $\Theta$ between the vertex set $C_v$ and the
vertex set $V(H)\setminus (W\cup C_v)$, so we can combine $\xi$ and
$\zeta$ into one proper $d$-coloring of the graph $\cal{G}$, also
the condition, that for any vertex $u\in A$ there are at least $q$
vertices of different colors in its neighborhood, holds for this
combined coloring. Thus we get a coloring of the graph $H$ we had
seeking for in the \textit{lemma \ref{l1}}, so we arrive at a
contradiction with assumption of the \textit{statement \ref{p8}}.

\bigskip

\textbf{2) $r\leq q^2-4$}

So from the vertex $v$ in the graph $\widetilde{I}$ it outcomes not
more than $q^2-4$ edges to the vertex set $R\cap V(\widetilde{I})$.
The degree of $v$ in the graph $\widetilde{I}$ is $d-q+q_1$. So from
the vertex $v$ it comes at least $d-q+q_1-r-l'$ edges to the set
$\{u_1,u_2,....,u_l\}$. Let us notice that if a vertex $u\in S$ in
the graph $\widetilde{\cal{G}}$ has an edge of multiplicity $k$,
then $d_{A}(u)\geq k-1$. We know that from $v$, the outcome is at
least $q^3+2q^2-q-8-q+q_1-r-l'$ edges to the vertices
$u_1,u_2,...,u_l$. Denote for all $i\in [1,l]$ as $d_i$ the
multiplicity of the edge $(v,u_i)$ in the graph $\widetilde{H}$.
Thus we know that $\sum\limits_{i=1}^{l}d_i\geq
q^3+2q^2-2q-8+q_1-r-l'$. Then we get
$$\sum\limits_{i=1}^{l}d_{A}(u_i)\geq
\sum\limits_{i=1}^{l}(d_i-1)\geq
q^3+2q^2-2q-8+q_1-r-l'-l~.\eqno(10)$$

By substituting inequality (5) and equality (6) into inequality (3)
we get

$$l'-q+q_1+\frac{1}{q+1}+\sum\limits_{i=1}^{l}d_{_{\widetilde{I}}}(u_i)-L(u_i)<0~.\eqno(11)$$

We know that  $d_{_{\widetilde{I}}}(u_i)-L(u_i)\geq d
-d_{\cal{G}}(u_i)-\frac{d_{A}(u_i)}{q+1}$ for all $i\in [1,l]$. By
applying inequality (1) we get $d
-d_{\cal{G}}(u_i)-\frac{d_{A}(u_i)}{q+1}\geq \lceil
\frac{d_{A}(u_i)}{q}\rceil -\frac{d_{A}(u_i)}{q+1}\geq
d_{A}(u_i)\frac{1}{q(q+1)}$. Thus
$$d_{_{\widetilde{I}}}(u_i)-L(u_i)\geq d_{A}(u_i)\frac{1}{q(q+1)}~.
\eqno(12)$$

Substitute (12) into (11) we get:

$$l'-q+q_1+\frac{1}{q+1}+\frac{1}{q(q+1)}\sum\limits_{i=1}^{l}d_{A}(u_i)<0~. \eqno(13)$$

Substitute (10) into (13):

$$l'-q+q_1+\frac{1}{q+1}+\frac{q^3+2q^2-2q-8+q_1-r-l'-l}{q(q+1)}<0~.$$

We know that $l'\geq 0$ and $q_1\geq 1$. Hence we have

$$l'(1-\frac{1}{q(q+1)})-(q-1)+\frac{1}{q+1}+\frac{q^3+2q^2-2q-8+1-r-l}{q(q+1)}<0~.$$

We also know that $r\leq q^2-4$ and by inequality (8) $l\leq q^2-3$.
Then

$$\frac{1}{q+1}+\frac{q^3+2q^2-2q-8+1-(q^2-4)-(q^2-3)-(q^3-q)}{q(q+1)}<0~.$$

So we get $\frac{1}{q+1}-\frac{q}{q(q+1)}<0$, i.e. we arrive at a
contradiction.

\medskip

Thus we proved inequality $(2')$ for the graph $\widetilde{I}$ and
so we proved the \textit{statement \ref{p8}}.

\end{proof}
\end{Proposition}

Return to the \textit{lemma \ref{l1}} proof. Let us begin for the
given sampling of sets $\lambda\in\Lambda$ to \textit{delete
recursively} vertices from $\widetilde{\cal{G}}$ while it is
possible. If graph $\widetilde{\cal{G}}$ is not a \textit{recursive}
one, then a graph $\widetilde{\cal{G}}(S'\cup R')$ will remain from
$\widetilde{\cal{G}}$, where $S'\subseteq S$, $R'\subseteq R$ and
$S'\neq \o$. Let us choose among all samplings from $\Lambda$ such a
sampling $\lambda$ that the value of $|S'|+|R'|$ achieves minimum on
it. Let us check up all assumptions of the \textit{statement
\ref{p8}} for sets $S'$ and $R'$. The unique non-trivial place in
this check is to verify inequality (2).

Since we can not \textit{delete recursively} from the graph
$\widetilde{H'}:=\widetilde{\cal{G}}(S'\cup R')$ any vertex, then
the degree of any vertex there is at least $d$. So by inequality (1)
$L(u)=d_{\cal{_G}}(u)+\frac{d_{_A}(u)}{q+1}\leq d$ for all $u\in B$,
at that $L(u)<d$ for all vertices $u\in S$. Since $S'\neq\o$ it is
clear that
$$\sum\limits_{u\in\widetilde{H'}}d_{\widetilde{H'}}(u)\geq
d|\widetilde{H'}|>\sum\limits_{u\in\widetilde{H'}}L(u).$$

So we will apply the \textit{statement \ref{p8}} to sets $S'$ and
$R'$, while we get such a vertex such, that it degree in
$\widetilde{H'}$ is less than $d$ (let us notice that we wan't do an
infinite number of \textit{regular non-degenerate changes} in regard
to the set $S'\cup R'$, since by any such a change we decrease the
value of $\sum\limits_{i=1}^{|A|}|S_i\cap S'|$). Due to the
\textit{statement \ref{p8}} we can as before to \textit{delete
recursively} all vertices from $B$ except $S'\cup R'$, and then we
can to delete recursively one extra vertex of degree less than $d$
from the set $S'\cup R'$. Thus we arrive at a contradiction with
minimality of $|S'|+|R'|$. Hence there is such a sampling of
\textit{permissible} sets that the graph $\widetilde{\cal{G}}$ would
be a \textit{recursive} one. Thus we proved the \textit{lemma
\ref{l1}} and finally we proved the \textit{theorem \ref{t1}} (see
\textit{remark \ref{r35}}).
\end{proof}

\end{Lem}

\end{proof}
\end{theorem}

\begin{flushleft}
\bf\large{Bibliography}
\end{flushleft}

[1] R.\,L.\,Brooks, \textit{On coloring the nodes of network}, Proc.
Cambridge Philos. Soc. 37 (1941) 194-197.

[2] L.\,Lovasz, \textit{On decomposition of graphs}, Studia Sci.
Math. Hungar. 1 (1966) 237-238.

[3] L.\,Hong-Jian, B.\,Montgomery, Hoifung Poon, \textit{Upper
Bounds of Dynamic Chromatic Number}, Ars. Combinatoria 68(2003), pp.
193-201.

[4] B.\,Montgomery, \textit{Dynamic Coloring}, Ph.D. Dissertation,
West Virginia University, 2001.

[5] J.\,A.\,Bondy, U.\,S.\,R.\,Murty, \textit{Graph Theory with
Applications}, American Elsevier, New York, 1976.

[6] X.\,Meng, L.\,Miao, Z.\, R.\,Li, B.\,Su, \textit{The Conditional
coloring numbers of Pseudo-Halin graphs}, Ars Combinatoria 79(2006),
pp. 3-9.

[7] Suohai Fan, Hong-Jian Lai, Jianliang Lin, Bruce Montgomery,
Zhishui Tao, \textit{Conditional colorings of graphs}, Discrete
Mathematic 16(2006), 1997-2004.
\end{document}